Pryjmak Mykola Volodymyrovych Ternopil Ivan Pul'uj National Technical University

# PERIODIC FUNCTIONS WITH THE VARIABLE PERIOD

The examples of rhythmical signals with variable period are considered. The definition of periodic function with the variable period is given as a model of such signals. The examples of such functions are given and their variable periods are written in the explicit form. The system of trigonometric functions with the variable period is considered and its orthogonality is proved. The generalized system of trigonometric functions with the variable period is also suggested; some conditions of its existence are considered.

**Introduction.** In applied research we sometimes deal with signals (phenomena, systems), which have rhythmical (oscillating, vibrating, cycling etc.) character. We understand the term «rhythmical signals» as signals for which strict periodicity (repetition) is absent, but their certain stochastic characteristics, such as expectation value and dispersion repeat periodically. The power system load and gas consumption diagrams appear to be rhythmic when they are examined on time intervals as compared to twenty-four hours. Incoming flows of most queuing systems and dataflow (traffic) in the Internet networks are referred to as rhythmic. The ambient temperature, illuminance change rhythmically with the period of 24 hours and in addition to the daily periodicity, the seasonal (annual) periodicity also occurs for the given meteorological factors. According to the classical approach "model-algorithm-program" which is recommended to use for the research of complex objects, the rhythmic signals study begins with their model substantiation. Having got the model, one can develop (or choose from the known) methods and algorithms of the signal processing. When rhythmic signals are determined or close to it, periodic functions can be used as their model, and methods can be used for their analysis . They are mostly based on the Fourier theory, amplitude and phase spectra of the periodic functions. There is a range of random processes and sequence to describe rhythmic signals taking into consideration their stochastic character. They are random processes periodically correlated [1] and sequences [2], periodic random processes (according to Slutsky) [2,3, p. 703-755], periodic white noises with continuous [4] and discrete [5] argument, markovian periodic processes and periodic Markov chains [6].

**Rhythmical signals with the variable period.** In addition to rhythmic signals whose period is considered to be constant, there are signals which can be referred to as rhythmic only symbolically. The patient's electrocardiogram is an excellent example. It was recorded just after physical activity, and it is analyzed for some period of time while pulse rate is coming to normal. The figures 1a-1c show three segments of the electrocardiogram, each of 3 seconds duration and recorded in certain time intervals after physical activity. The fig. 1a shows the electrocardiogram recorded in 60 seconds after physical activity, the fig. 1b and 1c recorded in 120 seconds and 180 seconds accordingly after physical exercises. As follows from the diagram analysis the electrocardiogram form repeats roughly both on every diagram and on different diagrams. It is evident though that the repetition period varies, namely increases, and eventually stabilizes, so that it corresponds to the reduction of pulse rate. The spirogram behavior must be similar to the cardiogram, which is also recorded after physical activity or other activator of the person's psychophysical state. The examples of similar signals can be given in the functioning of some engineering systems. It can be engine and diesel generators operation in transitive regimes, for example, after the external loading variations.

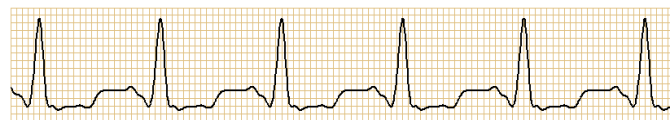

a) 1 minute after physical activity

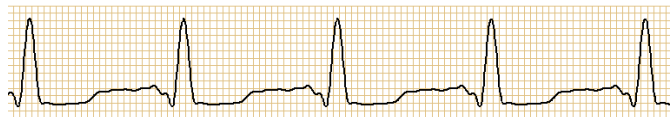

b) 2 minutes after physical activity

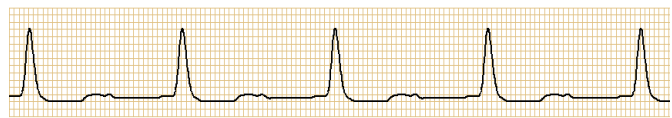

c) 3 minutes after physical activity

Fig. 1. The segments of the electrocardiogram recorded in different time intervals after physical activity.

The question arises, how shall we investigate rhythmic signals, processes, and phenomena with the variable period? While models and methods for periodic functions and processes have already been developed and became classical, only the first steps are made in research of rhythmic signals with variable period. The definition of periodic function with the variable period as well as the periodic random process with the variable period on its basis is given in [7]. The approaches to the digitization of such functions and processes are also described. However some questions concerning the properties of periodic functions with the variable period and the possibility of their studying by means of the Fourier methods remain urgent.

**The objective of the paper** is to give the definition of the periodic function with the variable period, examine the system of trigonometric functions with the variable period, and investigate the orthogonality of these functions, write down the form of trigonometric functions with the variable period in a general form and consider some conditions of their existence.

It should be noted that the function $g(x)$, $x \in I$, is called periodic if there is a number $T \neq 0$, when the values of $f(x)$ function from the domain $I$ are equal for all points $x$ and $x+T$, that is $g(x) = g(x+T)$. The number $T$ in this case is called the period of function. The function $I$ can coincide with the interval $(\infty, \infty)$, half-closed interval $[0, \infty)$, etc.

Let's return to the examples of the signals mentioned above. Their basic property is rhythm, and the period of rhythmicity also varies. The function with the variable period can be chosen as their model in the first approximation, the definition of which is given in [7].

**Definition 1.** The function $f(x)$, $x \in I$ is called periodic one with the variable period, if there is such function as $T(x) > 0$, when the values of $f(x)$ function from the domain $I$ are equal for all points $x$ and $x + T(x)$, that is

$$f(x) = f(x + T(x)). \tag{1}$$

The function $T(x)$ is called the variable period, and thus it should meet certain conditions, in particular it has to be continuous and its derivative $T'(x) > -1$.

If in (1) $T(x) = T = const$, we obtain the classical definition of the periodic function.

The example of the diagram of the period $T(x)$ is shown in figure 2. The functional period is equal to $T(x_1)$ in the point $x_1$ that is the value of the function in the points $x_1$ i $x_1 + T(x_1)$ is repeated: $f(x_1) = f(x_1 + T(x_1))$. In the point $x_2$ the number $T(x_2)$ is the period. From the given figure we can see that the function periods $f(x)$ are different in the points $x_1$ and $x_2$.

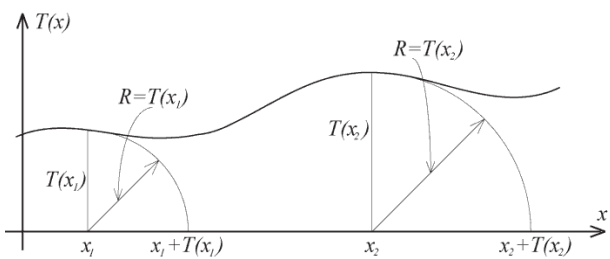

Fig. 2 The variable period $T(x)$, its value in points $x_1$ and $x_2$, and their corresponding points $x_1 + T(x_1)$ and $x_2 + T(x_2)$, in which the value of the function is repeated.

It should be noted that for the periodic function $g(x)$, $x \in I$ with the period $T$ the following equation is fulfilled

$$g(x) = g(x+T) = g(x-T). \tag{2}$$

For the periodic function $f(x)$ with the variable period the similar equation $f(x) = f(x+T(x)) = f(x-T(x))$ in general is not fulfilled and it is possible only in some instances. Therefore in order to find a point, placed in the period on the left with respect to $x$, in which the value of the function is repeated, the variable period $T^-(x)$ should be introduced, so that

$$f(x) = f(x - T^-(x)).$$

The relation for the variable periods $T(x)$ and $T^-(x)$ is fulfilled:

$$T(x) = T^-(x + T(x)), \tag{3a}$$

$$T^-(x) = T(x - T^-(x)). \tag{3b}$$

The function $T(x)$ as the variable period of the function $f(x)$ should meet certain conditions. According to the definition it should be continuous. It is also easy to show that its derivative should be greater than $-1$: $T'(x) > -1$. Indeed, as $x + \Delta x > x$, then for the corresponding points $x + T(x)$ and $x + \Delta x + T(x + \Delta x)$, taken in the period, the inequation $x + \Delta x + T(x + \Delta x) > x + T(x)$ or $T(x + \Delta x) - T(x) > -\Delta x$ should also be fulfilled. As follows from the last inequation $\dfrac{T(x+\Delta x) - T(x)}{\Delta x} > -1$. At the limiting transition we obtain that for the variable period $T(x)$ its derivative is

$$T'(x) > -1, \; x \in I. \tag{4}$$

From (4) it is obvious, that when the period $T(x)$ decreases on certain intervals, then this decrease should occur more slowly, than the function decrease $y(x) = -x$. Detection of the other properties of the variable period $T(t)$ requires separate research.

**Examples of the functions with the variable period and their variable periods.** Trigonometric functions $\sin x^\alpha, \cos x^\alpha, \alpha > 0, \alpha \neq 1$ are the simplest ones among the periodic functions with the variable period. Figure 3 shows the periodic function with the variable period $f_1(x) = \sin x^{3/4}, x \geq 0$, (diagram 1) and the function $f_2(x) = \sin x$ (diagram 2) to compare.

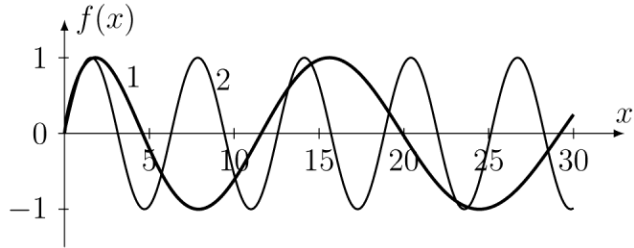

Fig.3. Function $f_1(x) = \sin x^{3/4}$ (diagram 1), $f_2(x) = \sin x$ (diagram 2).

As follows from the diagrams compared, the function period increases with the argument $x$ increase for the function $f(x) = \sin x^{3/4}$. Two periodic oscillations are located on the interval $[0, 30]$ for this function, whereas more than four oscillations are located on the same interval for the function $\sin x$.

For the sinusoidal functions with the variable period $\sin x^\alpha, \cos x^\alpha, \alpha > 0, \alpha \neq 1$ the following lemma is true.

**Lemma 1.** For the functions $\sin x^\alpha$ and $\cos x^\alpha$, $x \geq 0, \alpha > 0, \alpha \neq 1$, their variable periods being:

$$T_\alpha(x) = -x + \left(x^\alpha + 2\pi\right)^{1/\alpha}, \; x \geq 0, \quad (5a)$$

$$T_\alpha^- = x - \left(x^\alpha - 2\pi\right)^{1/\alpha}, \; x \geq (2\pi)^{1/\alpha}. \quad (5b)$$

These statements can be easily proved. Actually, the function $\sin x^\alpha$ repeats in the period $T_\alpha(x)$:

$$\sin\left(x + T_\alpha(x)\right)^\alpha = \sin\left(\left(x - x + \left(x^\alpha + 2\pi\right)^{1/\alpha}\right)^\alpha\right) = \sin\left(\left(x^\alpha + 2\pi\right)^{1/\alpha}\right)^\alpha = \sin\left(x^\alpha + 2\pi\right) = \sin x^\alpha,$$

$$\sin(x - T_\alpha^-(x))^\alpha = \sin\left(\left(x - x + (x^\alpha - 2\pi)\right)^{1/\alpha}\right)^\alpha = \sin\left(\left(x^\alpha - 2\pi\right)^{1/\alpha}\right)^\alpha = \sin(x^\alpha - 2\pi) = \sin x^\alpha.$$

Similar equations are true for the function $\cos x^\alpha$.

It should be noted that if there is no misunderstanding, an index $\alpha$ included in expressions $T_\alpha(x)$ can decrease.

Considering (5a) and (5b), for the function $\sin x^{3/4}$, $x \geq 0$, their variable periods are expressed by the formulas: $T(x) = -x + \left(x^{3/4} + 2\pi\right)^{4/3}$, $x \geq 0$; $T^-(x) = x - \left(x^{3/4} - 2\pi\right)^{4/3}$, $x \geq (2\pi)^{4/3} \approx 11.594$. The diagrams of these periods are shown in figure 4. The period $T = 2\pi$ is also shown for function $\sin x$ for comparison. As follows from the figure the period $T(x)$ increases with the increase of the argument. Having applied the formula (5a), we will find for example, that for the given $\alpha = 3/4$ the period gets the value $T(0) \approx 11.594$ in the point $x = 0$, for $x = 30$ the period $T(30) \approx 21.062$. As follows from the figure the period $T^-(x)$ decreases with the argument decreasing.

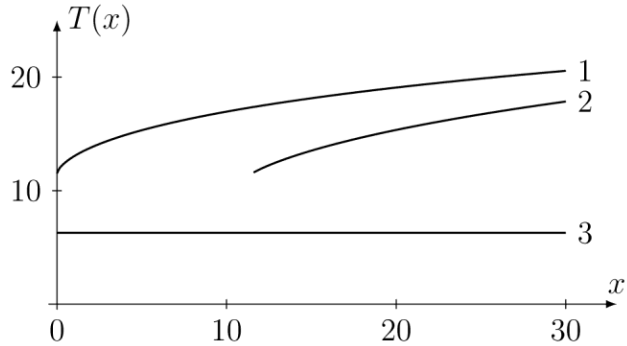

Fig. 4. The variable periods for the function $\sin x^{3/4}$, $x \geq 0$: $T(x) = -x + \left(x^{3/4} + 2\pi\right)^{4/3}$, $x \geq 0$ (diagram 1); $T^-(x) = x - \left(x^{3/4} - 2\pi\right)^{4/3}$, $x \geq (2\pi)^{4/3} \approx 11.594$ (diagram 2). The period $T = 2\pi$ for the function $\sin x$ ( diagram 3) is shown for the comparison.

Another example of the function with the variable period is shown in the figure 5. It is a diagram of the function $f_1(x) = \sin x^{4/3}$, $x \geq 0$. Provided that more than five and a half oscillations are located on the interval [0,15] of this sinusoid, and apparently, its period gradually decreases, then the function $f_2(x) = \sin x$ makes a bit more than two fluctuations.

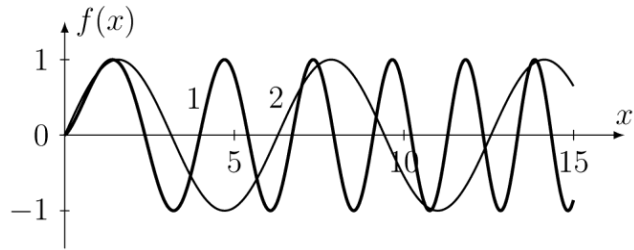

Fig.5. The diagrams of the function: $f_1(x) = \sin x^{4/3}$ (diagram 1); $f_2(x) = \sin x$ (diagram 2).

The diagrams of the variable periods for this function are shown in figure 6. Unlike the previous case the period $T(x)$ here decreases with the argument increasing. If, say, the period $T(0) = 3.968$ for $x = 0$, then the period $T(15) = 1.873$ in the point $x = 15$.

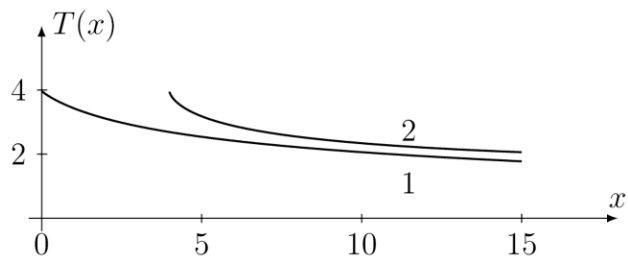

Fig. 6. The variable periods for the

function $f_1(x) = \sin x^{4/3}$, $x \geq 0$: $T(x) = -x + \left(x^{4/3} + 2\pi\right)^{3/4}$, $x \geq 0$ (diagram 1); $T^- = x - \left(x^{4/3} - 2\pi\right)^{3/4}$,

$$x \geq (2\pi)^{3/4} \approx 3.968 \text{ (diagram 2)}.$$

**Orthogonality of the trigonometric system of functions with the variable period.** Similar to trigonometric system of the functions 1, 1, $\sin mx$, $\cos mx$, $m = 1, 2, ...$, which is orthogonal at any segment of the $2\pi$ length, the natural question arises concerning the existence of the orthogonal systems of functions with the variable period, in particular in the set of trigonometric functions with the variable period. The following theorem appears to be true.

**The theorem 1.** The system of trigonometric functions with the variable period

$$\sin mx^\alpha, \cos mx^\alpha, x \geq 0, \alpha > 0, m = 1, 2, ..., \tag{6}$$

is orthogonal in any segment $[x_0, x_0 + T_\alpha(x_0)]$, $x_0 \in I$, with the weight function $\rho(x) = x^{\alpha-1}$.

**The proof.** Orthogonality can be easily proved by direct calculation, for example, for the functions $\sin mx^\alpha$, $m = 1, 2, ...$, and $\cos nx^\alpha$, $n = 1, 2, ...$, integral

$$\int_{x_0}^{x_0 + T_\alpha(x_0)} x^{\alpha-1} \sin mx^\alpha \cos nx^\alpha dx = \frac{1}{\alpha} \int_{x_0}^{x_0 + T_\alpha(x_0)} \sin mx^\alpha \cos nx^\alpha dx^\alpha =$$

$$\frac{1}{2\alpha} \int_{x_0}^{x_0 + T_\alpha(x_0)} \left( \sin(m+n)x^\alpha + \sin(m-n)x^\alpha \right) dx^\alpha =$$

$$= -\frac{1}{2\alpha} \left[ \frac{1}{m+n} \cos(m+n)x^\alpha + \frac{1}{m-n} \cos(m-n)x^\alpha \right]_{x_0}^{x_0+T_\alpha(x_0)} =$$

$$= -\frac{1}{2\alpha} \left\{ \frac{1}{m+n} \left[ \cos(m+n)(x_0 + T_\alpha(x_0))^\alpha - \cos(m+n)x_0^\alpha \right] + \right.$$

$$\left. + \frac{1}{m-n} \left[ \cos(m-n)(x_0 + T_\alpha(x_0))^\alpha - \cos(m-n)x_0^\alpha \right] \right\}.$$

Since according to (5a) the period $T_\alpha(x) = -x + (x^\alpha + 2\pi)^{1/\alpha}$, then $[x_0 + T_\alpha(x_0)]^\alpha = \left[ x_0 - x_0 + (x_0^\alpha + 2\pi)^{1/\alpha} \right]^\alpha = x_0^\alpha + 2\pi$. Taking it into consideration in the last expression we have the following

$$\int_{x_0}^{x_0 + T_\alpha(x_0)} x^{\alpha-1} \sin mx^\alpha \cos nx^\alpha dx = -\frac{1}{2\alpha} \left\{ \frac{1}{m+n} \left[ \cos(m+n)(x_0^\alpha + 2\pi) - \cos(m+n)x_0^\alpha \right] + \right.$$

$$\left. + \frac{1}{m-n} \left[ \cos(m-n)(x_0^\alpha + 2\pi) - \cos(m-n)x_0^\alpha \right] \right\} =$$

$$= -\frac{1}{2\alpha} \left\{ \frac{1}{m+n} \left[ \cos\left[(m+n)x^\alpha + 2\pi(m+n)\right] - \cos(m+n)x^\alpha \right] + \right.$$

$$\left. + \frac{1}{m-n} \left[ \cos\left[(m-n)x^\alpha + 2\pi(m-n)\right] - \cos(m-n)x^\alpha \right] \right\} =$$

$$= -\frac{1}{2\alpha} \left\{ \frac{1}{m+n} \left[ \cos(m+n)x^\alpha - \cos(m+n)x^\alpha \right] + \frac{1}{m-n} \left[ \cos(m-n)x^\alpha - \cos(m-n)x^\alpha \right] \right\} = 0.$$

We can check the orthogonality of the functions $\sin mx^\alpha$ and $\sin nx^\alpha$, the functions $\cos mx^\alpha$ and $\cos nx^\alpha$ in the same way provided that $m \neq n$. If $m = n$, then it is easy to prove that

$$\int_{x_0}^{x_0+T(x_0)} x^{\alpha-1} \sin^2 mx^\alpha dx = \frac{\pi}{\alpha}, \quad \int_{x_0}^{x_0+T(x_0)} x^{\alpha-1} \cos^2 mx^\alpha dx = \frac{\pi}{\alpha}. \tag{7}$$

Taking into account the standard specifications, the expression (7) can be given as

$$\|\sin mx^\alpha\|^2 = \frac{\pi}{\alpha}, \quad \|\cos mx^\alpha\|^2 = \frac{\pi}{\alpha}, \quad x \geq 0, m = 1, 2, \ldots. \tag{7a}$$

The theorem is proved.

From the results given above some preliminary conclusions can be made. Every $\alpha > 0$ can match to the orthogonal system of the trigonometric functions

$$\sin mx^\alpha, \cos mx^\alpha, \quad x \geq 0, \alpha > 0, m = 1, 2, \ldots, \tag{6}$$

with the variable period

$$T_\alpha(x) = -x + \left(x^\alpha + 2\pi\right)^{1/\alpha}, \quad x \geq 0. \tag{5a}$$

Meanwhile the segment of the orthogonality is any segment $[x, T_\alpha(x)]$ at length $T_\alpha(x)$, and the standard of the functions is defined according to the formulas

$$\|\sin mx^\alpha\|^2 = \frac{\pi}{\alpha}, \quad \|\cos mx^\alpha\|^2 = \frac{\pi}{\alpha}, \quad x \geq 0, m = 1, 2, \ldots. \tag{7a}$$

If we take $\alpha = 1$ for the system of the trigonometric function with the variable period (6), then it changes into common system of trigonometric functions $\sin mx, \cos mx, m = 1, 2, \ldots$ As results from (5a) the period $T = 2\pi$, and from (7a) it is clear, that the norms $\|\sin mx\|^2 = \pi$, $\|\cos mx\|^2 = \pi$. The last comments show, that trigonometric system of functions $\sin mx, \cos mx, m = 1, 2, \ldots$ can be considered as a partial case of the trigonometric system of the functions (6) with the variable period (5a) provided that $\alpha = 1$.

**The generalized system of trigonometric functions with the variable period.** The system of functions can be referred to as the generalisation of the system of trigonometric functions (6)

$$\sin mg(x), \cos mg(x), \quad x \in I, m = 1, 2, \ldots, \tag{8}$$

$g(x)$ being – continuous function, which should meet certain conditions. In order to specify some of these conditions, at first we will consider the problem of the variable period $T_g(x)$, corresponding to the $g(x)$ function. The period $T_g(x)$ is considered to provide the fulfillment of the equation

$$\sin\left[g\left(x+T_g(x)\right)\right]=\sin\left[g(x)+2\pi\right],$$

that is,

$$g\left[x+T_g(x)\right]=g(x)+2\pi. \tag{9}$$

Suppose there is inverse function $g^{-1}(y)$ to the function $g(x)$. Having applied this function (9), we obtain $x+T_g(x)=g^{-1}\left[g(x)+2\pi\right]$. Thus the variable period is calculated according to the formula

$$T_g(x)=-x+g^{-1}\left[g(x)+2\pi\right]. \tag{10a}$$

The expression (4) shows that the derivative of the variable period should be greater than $-1$. It is also true for the variable period $T_g(x)$, that is, according to (10), the following inequation is:

$$T'_g(x)=-1+\left(g^{-1}\left[g(x)+2\pi\right]\right)' >-1. \tag{11}$$

As follows from the last inequation the function $g(x)$ and its inverse function $g^{-1}(\circ)$ should provide the fulfillment of derivative of the variable period $T_g(x)$, which is expressed through these functions (11).

Similar to (9) for the variable period $T_g^-(x)$ the equation $g\left[x-T_g^-(x)\right]=g(x)-2\pi$ should be fulfilled. Having applied the inverse function to it, we obtain

$$T_g^-(x)=x-g^{-1}\left[g(x)-2\pi\right]. \tag{10b}$$

As the example, the figure 7 shows the diagram of the function $f_1(x)=\sin g(x)$ for $g(x)=x^{4/3}+1.2\sin x$, that is the diagram of the function $f_1(x)=\sin\left(x^{4/3}+1.2\sin x\right)$ and the diagram of the function $f_2(x)=\sin x$ for comparison.

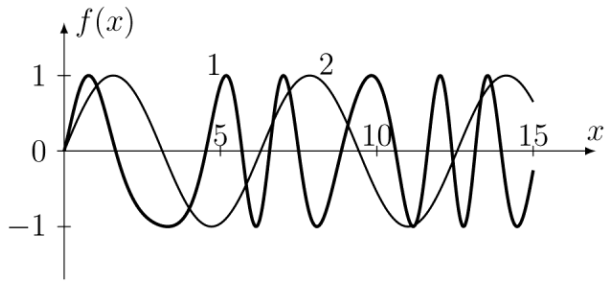

Fig. 7. The function $f_1(x)=\sin\left(x^{4/3}+1.2\sin x\right)$ (diagram 1), $f_2(x)=\sin x$ (diagram 2).

As follows from the behaviour of the function diagram $f_1(x) = \sin\left(x^{4/3} + 1.2\sin x\right)$ its compression (at the expense of period reduction) occurs to certain "indignation", which is caused by the influence of the sinusoid $1.2\sin x$ on the function $x^{4/3}$ behavior.

It should be mentioned, that when the function is $f(x) = x$ in the system (8), we obtain the common system of trigonometric functions. Provided $f(x) = x^\alpha, \alpha > 0$, we obtain the trigonometric system of functions (6) with the variable period (5a).

To sum up we should note, that in the further research of the trigonometric system of functions (8), some additional conditions will be revealed in addition to the conditions, mentioned above which should fulfill the functions $g(x)$ and $T_g(x)$. It is necessary to indicate, that in addition the problem considered above dealing with the finding of expressions for the variable periods $T_g(x)$ and $T_f^-(x)$ (formulas (10a) and (10b) accordingly) by means of the function $g(x)$, the problem of finding the function $g(x)$ is very important for the applied problems in particular, provided the variable period $T_g(x)$ is known.

**Conclusions.** The given definition of the periodic function with the variable period and the examined system of the orthogonal trigonometric functions with the variable period expand considerably the concept and meaning of the periodic functions and trigonometric systems with the constant period. The examined functions appear to be not only of the theoretical interest, but also they can be used in applied research of rhythmic signals and phenomena with the variable period. The research of the functions of new kind is apparently waiting for its continuation.